\documentclass[a4paper]{article} 
\usepackage{a4wide}
\usepackage{amsmath}
\usepackage{amssymb}
\usepackage{comment}
\usepackage{hyphenat}
\usepackage{multicol}
\usepackage{tabularx}

\usepackage{subfigure} 
\usepackage{caption}
\usepackage{graphicx}

\usepackage{hyperref}
\usepackage{bm}
\usepackage{float}
\usepackage{xspace}
\usepackage{fancyvrb}

\newcommand{\LRb}[1]{\left(#1\right)}
\newcommand{\LRB}[1]{\Bigl(#1\Bigr)}

\newcommand{\LRs}[1]{\left[#1\right]}    

\newcommand{\SET}[1]{\left\{#1\right\}}

\newcommand{\FRAC}[2]{\frac{\displaystyle{#1}}{\displaystyle{#2}}}
\newcommand{\MOD}[1]{\left|#1\right|}

\newcommand{\NORM}[1]{\left\|#1\right\|}
\newcommand{\LE}{{\leqslant}}
\newcommand{\GE}{{\geqslant}}
\newcommand{\EQ}{{=}}
\newcommand{\eps}{\varepsilon}

\newcommand{\eqdef}{\doteq}

\newcommand{\RR}{\mathbb{R}}

\newcommand{\IN}{{\in}}

\newcommand{\quot}[1]{``{#1}''}

\newcommand{\SVF}{{SvF}\xspace}

\DeclareMathOperator*{\Argmin}{\mathsf{Argmin}}

\providecommand{\keywords}[1]{\textbf{Keywords:}~#1}
\date{\today ~(revised)}
  
\begin{document}
 
%
%
%
%
%
%

\title{Balanced Identification as an Intersection of Optimization and Distributed Computing}
\author{Alexander V. Sokolov\footnote{Center for Distributed Computing, Institute for Information Transmission Problems (Kharkevich institute),  Bolshoy Karetny per. 19-1, Moscow 127051, Russia, \texttt{alexander.v.sokolov@gmail.com, vv\_voloshinov@iitp.ru}} \and Vladimir V. Voloshinov\footnotemark[\value{footnote}]}



  \begin{abstract}
{
Technology of formal quantitative estimation of the conformity of the mathematical models to the available dataset is presented. 
Main purpose of the technology is to make easier the model selection decision-making process for the researcher.  
The technology is a combination of approaches from the areas of data analysis, optimization and distributed computing including: cross-validation and regularization methods, algebraic modeling in optimization and methods of optimization, automatic discretization of differential and integral equation, optimization REST-services. 
The technology is illustrated by a demo case study. General mathematical formulation of the method is presented. It is followed by description of the main aspects of algorithmic and software implementation. 
Success story list of the presented approach already is rather long. Nevertheless, domain of applicability and important unresolved issues are discussed.}
\end{abstract}

%
%
  
\maketitle
\section{Introduction}

The same object, process or phenomenon can be described by different
mathematical models. The model selection among possible candidates remains one of the  important problems in applied mathematics.
A complexity, dimension and detail level of  mathematical description are defined by both the qualitative knowledge of the object (which can be formalized as mathematical statements), and the availability of quantitative data, their volume, detail, reliability and accuracy. 
The more complex is a phenomenon under study and the corresponding mathematical model, the more detailed and reliable the measurement should be.
Usually to find a balance between a model complexity and an available measurement 
an easy approach is used: various models and datasets are tested until the built model meets the goals of research.
There is no universally accepted approach here still. 
In each case, various parametrization, different identification and verification procedures  and variety of software are used. 
In addition, the model selection criterion is not formalized and the selection of the particular model is often made subjectively.

What we need are: a technology that simplifies the model selection procedure for the researcher; 
a unified method that gives a quantitative estimation of the conformity of the mathematical description to the available dataset for each model under study.

This paper presents a method of a \textit{balanced identification} (also referred to as \SVF-method, \textit{simplicity vs fitting}). The method has been developed to fulfil above requirements.
Besides formal description a software implementation of the \SVF-method is presented. This implementation is called the \SVF-technology. Regardless that the method is rather mature and was successfully used in a number of various research for a few years, its general mathematical scheme has been published recently \cite{bib:sokolov2018choice}. The basis for the \SVF-method are nonparametric identification with \textit{regularization} \cite{bib:tikhonov1980on} and well-known \textit{cross-validation} \cite{bib:kohavi1995study,bib:hastie2009elements,bib:kuhn2013applied}. Regularization enables to search the desired balance between simplicity and fitting to the dataset. Cross-validation error is used as a quantitative estimation of the model compliance to the dataset.

In general, \SVF-method may be applied to both parametric and nonparametric identification (where the model under study contains unknown functions). 
Depending on the class of the models, 
the variational problems arising here can be interpreted
as problems of variational spline approximation \cite{bib:rozhenko2005theory}, 
nonparametric regression \cite{bib:hardle1990applied},
predictive modeling \cite{bib:kuhn2013applied}, 
machine learning \cite{bib:hastie2009elements}, and others.

There is a trend in the research on nonlinear regression to create algorithms of \quot{automatic} composition of the best regression function from a predefined set of basic functions \cite{bib:strijov2010nonlinear}. Here, the object under study is considered as a black-box and the best regression remains a purely phenomenological description. \SVF-technology is not a fully automated routine. It enables quantitative comparison of a given set of models, but it is the researcher who make modifications of the model (e.g. include/exclude model constraint equations). Moreover, \SVF-method can be applied for structural mathematical models describing internal structure of the phenomenon, e.g. via algebraic or integro-differential equations. One can say that \SVF-method generalizes the concept of regression with additional constraints. Known example of that is a monotonic regression \cite{bib:sysoev2019smoothed}, where a monotonic response function should be found.

For a given model the \SVF-method requires solving a bilevel optimization problem, where at lower-level, we have a number of independent mathematical programming problems to get the value of objective function of the upper-level problem. Finaly, the optimal values of cross-validation error and regularization coefficients are obtained for a given model. Root-mean-square error will be presented for reference only. 
Bilevel optimization is a well-known hard challenge for numerical methods \cite{bib:dempe2002bilevel}. Current implementation of \SVF-method is based on surrogate optimization and on the ideas of \textit{Pshenichny-Danilin linearisation method} \cite{bib:pshenichnyi1993linearization}.

Because at the lower-level we have a number of independent problems, it is possible to increase performance of the technology by distributed computing. 
One should mention that optimization in distributed computing environment is already used in regression analysis, e.g. \cite{bib:boyd2011distributed}. 
In \SVF-technology the basic tool for solving independent mathematical programming problems in parallel is the Everest optimization service \cite{bib:smirnov2016distributed,bib:smirnov2018dd} based on Everest toolkit \cite{bib:EverestISPDC2015}, \mbox{\url{http://everest.distcomp.org}}.

The paper is organized as follows. The next section demonstrates an example of applying the \SVF-technology to the modeling of an oscillator with friction. 
Section \ref{sec:svf_outline} contains mathematical formulation of the basic concepts of \SVF-method including the formulation of the main bilevel optimization problem.
Section \ref{sec:tech} presents main aspects of the \SVF--method software implementation: usage of Pyomo package (Python Optimization Modeling Objects), \mbox{\url{http://www.pyomo.org}}; special symbolic notation that simplifies formulation of the model and regularization rule; discretization of differential and integral equations if they are present in the model; usage of Everest optimization service. 
Successful use of \SVF-technology in various researches are presented in Section \ref{sec:use_cases}.
Some important open issues related the the current implementation of \SVF-technology are discussed in Section \ref{sec:concl} followed by Acknowledgments. 

\newcommand{\Dx}{{\Delta}{x}}
\newcommand{\kmax}{{k_{\max}}}
\newcommand{\tmax}{{t_{\max}}}
\newcommand{\tmin}{{t_{\min}}}
\newcommand{\xmax}{{x_{\mathsf{1}}}}
\newcommand{\xmin}{{x_{\mathsf{0}}}}
\newcommand{\xaKk}[1]{{x^{#1}_{K{\setminus}k}}}
\newcommand{\faKk}[1]{{f^{#1}_{K{\setminus}k}}}
\newcommand{\sigmA}{{\sigma({\alpha})}}
\newcommand{\rmse}{\mathrm{rmse}}

\section{Demonstrative example}\label{sec:oscill}
Before proceeding with formal description of the method it is worth to consider a demonstrative \quot{use case} on the example of a classical damped oscillator described by the following function of time $t$:
\begin{equation}\label{eq:oscill}
x(t)=\sin\LRb{\frac{\sqrt{4k{-}\mu}}{2}t}{\cdot}\exp\LRb{-\frac{\mu}{2}t} + \Dx,
\end{equation}
\noindent where: $x(t)$ is a trajectory of oscillation; 
$k\EQ1.56$ (elastic coefficient); $\mu\EQ0.4$ (friction factor); $\Dx\EQ1.2$ (initial displacement of oscillator)

\begin{figure}[h]
\includegraphics[scale=0.43]{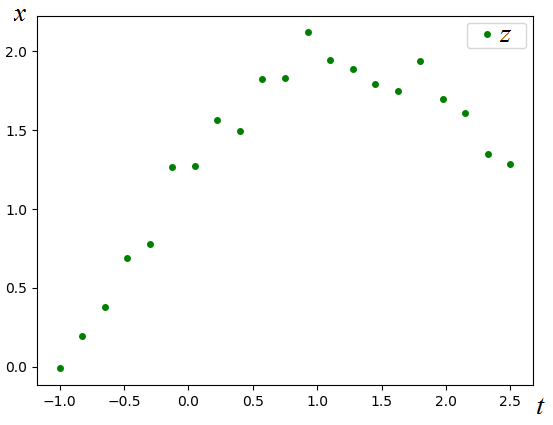}
\caption{Raw data of demo example}
\label{fig:oscill_D}
\end{figure}

Assume that we do not know in advance expression (\ref{eq:oscill}). 
What we have is a dataset (measurement series) on time interval $T$ (here and below $0{:}n$ denotes the set of numbers 
$\SET{0,1,2,{\dots},n}$):
\begin{equation}\label{eq:oscill_D}
\begin{array}{l}
D\EQ\SET{(z_k,t_k): k\IN K}, K\EQ0{:}\kmax,\\
t_k\EQ\tmin{+}k h_t, ~h_t\EQ\frac{\tmax{-}\tmin}{\kmax}\\
\tmin\EQ-1.0,~\tmax\EQ2.5,~z_k\EQ x(t_k){+}\eps_k,
\end{array}
\end{equation}
\noindent $\eps_k$ is a random error with zero mean and variance 0.1. The data distorted in this way are presented in Fig. \ref{fig:oscill_D}.
Hereinafter the following notation for intervals of possible values of $t$ and $x(t)$ will be used ($X$ is chosen with some margin):
\begin{equation}\label{eq:oscill_TX}
\begin{array}{l}
\!T\EQ\LRs{\tmin,\tmax}, X\!\EQ\LRs{{{-}}0{.}1,2{.}2}{\supset}\!\LRs{\min\limits_{k\IN K}{z_k},\max\limits_{k\IN K}{z_k}}.
\end{array}
\end{equation}

Let's set a problem of determining the equation of motion (\ref{eq:oscill}) by the method of balanced identification via given measurements (\ref{eq:oscill_D}). For that, consider several mathematical models and obtain corresponding CV-error, which will be used as a criterion when comparing models.

\begin{figure}[h]
\centering
\subfigure[$\alpha{\to}{+}\infty$]{{\includegraphics[width=0.31\linewidth]{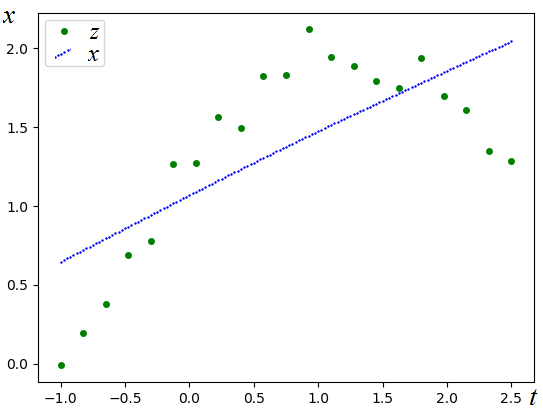}}} 
\subfigure[$\alpha{\to}0$]{{\includegraphics[width=0.31\linewidth]{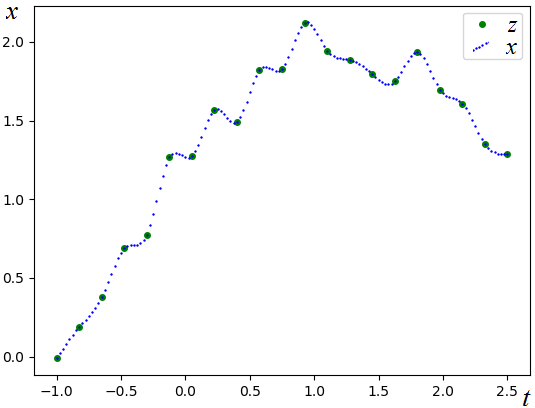}}} 
\subfigure[$\alpha-\mbox{fixed}$]{{\includegraphics[width=0.31\linewidth]{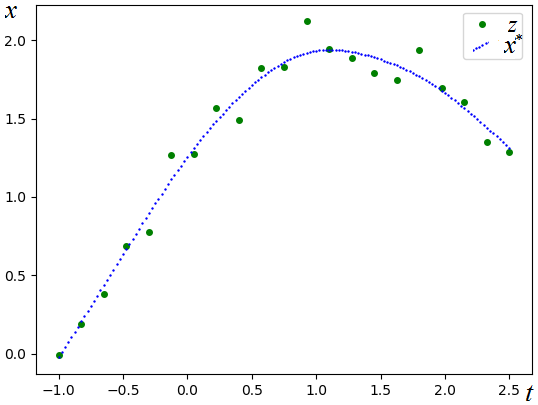}}}
\caption{Three variants of data approximation by a function}\label{fig:oscill_m0}
\end{figure}

One should say that all subsequent models of oscillator produce variational optimization problems. Direct solving of such problems is challenging because off-the-shelf optimization solvers (at least those, which are stable and \quot{mature}) can not handle integral and differential equations. To get a numerical solution all these problems are replaced with finite dimensional mathematical programming problems. To do that the discretization of differential and integral equations is applied (see details in Section \ref{sec:tech}).
\paragraph*{Model 0. Spline approximation of the function $x(t)$.} Take the simplest model as a twice differentiable function:
\begin{equation}
x(\cdot)\IN \mathrm{C}^2(T),
\end{equation}
\noindent and consider the following optimization problem (w.r.t. variable $x(\cdot)$), which depends on the regularization coefficient $\alpha\GE0$
\begin{equation}\label{eq:oscill_svf_m0}
\begin{array}{l}
\!\!\!\!\!\!F(x,K,\alpha){\eqdef}\frac{1}{\MOD{K}}\!\!\sum\limits_{k\IN K}\!\!\LRb{z_k{-}x(t_k)}^2\!{+}
\alpha\!\int\limits_{T}\!\!\LRb{\frac{d^2x}{dt^2}}^2\!\!\!dt\to\!\!\!\!\min\limits_{x{\in} \mathrm{C}^2(T)}. 
\end{array}
\end{equation}
\noindent For any fixed $\alpha$ (see \cite{bib:rozhenko2005theory,bib:hastie2009elements}) the problem (\ref{eq:oscill_svf_m0}) has unique solution, which is a cubic spline for the given set of points. It is presented in Fig. \ref{fig:oscill_m0} (c) and looks rather reasonable. 

Figure (a) (linear function $x(t)$) corresponds to the case $\alpha{\to}{+}\infty$, when the 2nd regularization term of objective function (\ref{eq:oscill_svf_m0}) suppresses the 1st one and the optimal solution tends to be a linear regression given by \textit{least squares method} (\textbf{too rough approximation}). In the case of (b) $\alpha{\to}0$ we have a problem of spline interpolation: find a curve with minimal integral curvature passing exactly through given points (\textbf{obvious overfitting}).  

It is required to find the \quot{best} value of $\alpha$ at which the \quot{smooth} model function passes \quot{close enough} to the measurements, {\bf giving optimal balanced approximation} and smoothing out random errors (as in Fig. \ref{fig:oscill_m0}, (c)). The choice of an optimal coefficient $\alpha$ can be made by minimizing the value of cross-validation error \cite{bib:kohavi1995study,bib:hastie2009elements,bib:kuhn2013applied}.

The simplest variant of the cross-validation (CV) procedure (\textit{leave-one-out}) is to use a one-point test sample (with number $k$) and a training sample consisting of the remaining set of measurements $K{\setminus}k$ (set (\ref{eq:oscill_D}) with one removed element $(z_k, t_k)$). 
Formally, for a fixed $\alpha$  and every $k$ (see (\ref{eq:oscill_xaKk})): find a solution of the following optimization problems, calculate approximation error on the test samples and get the average CV-error (which depends on $\alpha$):
\begin{equation}\label{eq:oscill_xaKk}
\begin{array}{l}
\xaKk{\alpha}{\eqdef}\Argmin\limits_{x(\cdot)}\SET{F\LRb{x(\cdot),K{\setminus}k,\alpha}}, k\IN K,\\
\sigmA\eqdef\sqrt{\FRAC{1}{\MOD{K}}\sum\limits_{k\IN K}\LRb{z_k{-}\xaKk{\alpha}(t_k)}^2}.
\end{array}
\end{equation}

%

\begin{figure}[H]
\centering
\includegraphics[width=0.6\linewidth]{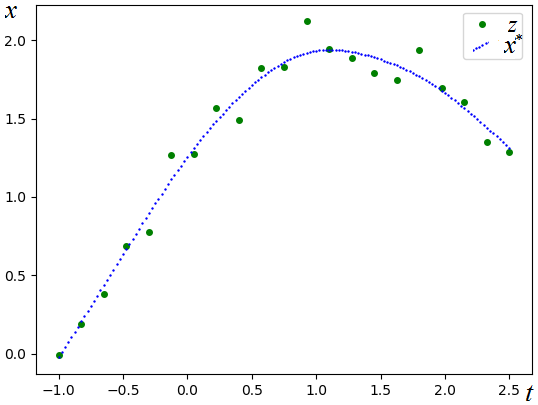}
\caption{Model 0. \SVF-approximation as cubic spline approximation}\label{fig:oscill_m0_svf}
\end{figure}

The best $\alpha^*$ minimises the value of CV-error: 
\begin{equation}\label{eq:oscill_aOpt}
\begin{array}{l}
\alpha^{*}\eqdef\Argmin\limits_{\alpha\GE0}\sigmA. 
\end{array}
\end{equation}
\noindent So, $\alpha^*$ is a solution of bilevel optimization problem: (\ref{eq:oscill_aOpt}) at upper-level 
and $\MOD{K}$ independent problems (\ref{eq:oscill_xaKk}) at lower-level. 
The value $\sigma^*{\eqdef}\sigma(\alpha^*)$ is called \textit{modeling CV-error}. 
It will be used as a quantitative measure of model's correspondence to the available experimental data:
\begin{equation}\label{eq:oscill_CVerror}
\begin{array}{l}
\sigma^{*}\eqdef\sqrt{\FRAC{1}{\MOD{K}}\sum\limits_{k\IN K}\LRb{z_k{-}\xaKk{\alpha^{*}}(t_k)}^2}.
\end{array}
\end{equation}
Finally, the $\alpha^*$ gives the sought model function $x^*(\cdot)$ and its \textit{root-mean-square error} (for reference only):
\begin{equation}\label{eq:oscill_xOpt}
\begin{array}{l}
x^*\eqdef x^{\alpha^*}_K\eqdef\Argmin\limits_{x(\cdot)}\SET{F\LRb{x(\cdot),K,\alpha^*}},\\
\rmse^*=\sqrt{\FRAC{1}{\MOD{K}}\sum\limits_{k\IN K}\LRb{z_k{-}x^*(t_k)}^2}.
\end{array}
\end{equation}
\noindent Function obtained for Model 0 is shown in Fig. \ref{fig:oscill_m0_svf}
(the same as in Fig. \ref{fig:oscill_m0}c)

Cubic spline approximation is one of the form of nonparametric regression, but other nonparametric models may be considered as well. 
Because the object under study seems to be a dynamical one, the following models include differential equation with unknown characteristics.

\paragraph*{Model 1. 1st order differential equation.} Apply \SVF-method for the model with two unknown functions $x(t)$, $f(t)$ related by the following differential equation (see definition of $T,X$ in (\ref{eq:oscill_TX})):
\begin{equation}\label{eq:oscill_m1}
\frac{dx}{dt}\EQ f\LRb{x(t)}, ~x(\cdot)\IN\mathrm{C}^2(T), ~f(\cdot)\IN \mathrm{C}^2(X), ~t\IN T.
\end{equation}

\begin{figure}[H]
\centering
\includegraphics[width=0.7\linewidth]{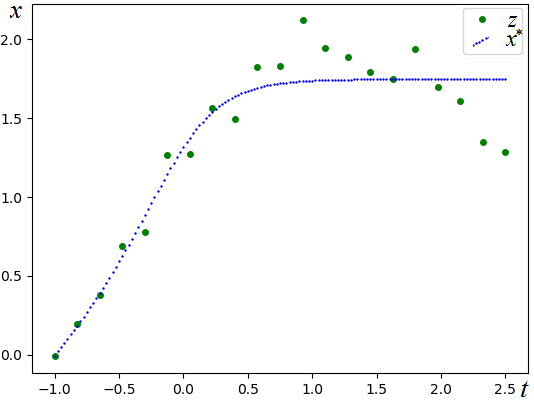}
\caption{Model 1. 1st order ODE: optimal $x^*(t)$ and raw data.}\label{fig:oscill_DE1_m1}
\end{figure}

In this case, it is reasonable to characterize complexity of the model by smoothness of function $f(x)$, responsible for dynamics, not $x(t)$. The main criterion (with regularization term) will be the following (hereinafter we'll use the same notation $F(\cdot)$, but with another argument list):
\begin{equation}\label{eq:oscill_F_m1}
\begin{array}{l}
\!\!\!\!F(x,f,K,\alpha){\eqdef}\frac{1}{\MOD{K}}\!\!\sum\limits_{k\IN K}\!\!\LRb{z_k{-}x(t_k)}^2\!{+}
\alpha\int\limits_{X}\!\!\LRb{\frac{d^2f(x)}{dx^2}}^2\!\!dx. 
\end{array}
\end{equation}
\noindent The \quot{balanced} solution $\LRb{x^*,f^*}$ for the dataset (\ref{eq:oscill_D}) corresponds to the $\alpha^*$ that minimizes the value of CV-error: 
\begin{equation}\label{eq:oscill_xfaKk}
\hspace{-3.75mm}\begin{array}{l}
\LRb{\xaKk{\alpha},\faKk{\alpha}}\!{\eqdef}\!\Argmin\limits_{x(\cdot),f(\cdot)}\SET{\!F\!\!\LRb{x,\!f,\!K{\setminus}k,\!\alpha}{:} \frac{dx}{dt}\EQ f\!(x(t))}\!\!,\\
\sigmA\eqdef\sqrt{\FRAC{1}{\MOD{K}}\sum\limits_{k\IN K}\LRb{z_k{-}\xaKk{\alpha}(t_k)}^2},\\
\alpha^{*}\eqdef\Argmin\limits_{\alpha\GE0}\sigmA,
\end{array}
\end{equation}

\noindent and $x^*,f^*$ are determined as following:
\begin{equation}\label{eq:oscill_xfOpt}
\begin{array}{l}
\LRb{x^*,f^*}{\eqdef}\Argmin\limits_{x(\cdot),f(\cdot)}\SET{F\LRb{x,f,K,\alpha^*}{:}\frac{dx}{dt}\EQ f\!\LRb{x(t)}}.\\
\end{array}
\end{equation}
\noindent The root-mean-square error is calculated as in (\ref{eq:oscill_xOpt}). See values of $\sigma^*$ and $\rmse^*$ in the 2nd row of Table \ref{tab:oscill} and the plot of function $x^*(t)$ in Fig. \ref{fig:oscill_DE1_m1}. All of them are obviously worse than those given by the Model 0.

Replacement of Model 0 with Model 1 (based on a solution of ODE (\ref{eq:oscill_m1}) instead of an arbitrary smooth function)
reduces the possibility of describing the data (\ref{eq:oscill_D}) by the function $x(t)$. This is so, because the set of feasible solutions of equation (\ref{eq:oscill_m1}) is too narrow.
It is easy to see that only a monotonic function $x(t)$ may be a solution of (\ref{eq:oscill_m1}). It is obviously insufficient for approximation of raw data in Fig. \ref{fig:oscill_D}. It is the reason why the CV-error $\sigma^*$ and $\rmse^*$ in the 2nd row of Table \ref{tab:oscill} are bigger than those in the 1st row. We conclude that Model 1 is worse than Model 0.

\paragraph*{Model 2. 2nd order differential equation.} 
Let's take a model with three unknown functions $x(t)$, $v(t)$ and $f(x,v)$ related by the following 2nd order ODE:
\begin{equation}\label{eq:oscill_m2}
\begin{array}{l}
\frac{d^2x}{dt^2}\EQ f\LRB{x(t),v(t)}, ~\frac{dx}{dt}\EQ v(t), V{\eqdef}\LRs{{-}1{.},1{.}5},\\
x(\cdot),v(\cdot)\in\mathrm{C}^2(T), f(\cdot)\IN \mathrm{C}^2(G){\subset}\RR^2,~G\EQ X{\times}V.\\
\end{array}
\end{equation}
\noindent The interval $V$ is wide enough to hold all possible values of $\dot{x}(t)$.

Now, the regularization term of \SVF-criterion includes two coefficients $\alpha_x$, $\alpha_v$ and partial derivatives of $f(x,v)$:
\begin{equation}\label{eq:oscill_Fxf_m2}
\begin{array}{l}
\!\!F(x,f,K,\alpha_x,\alpha_v){\eqdef}\frac{1}{\MOD{K}}\!\!\sum\limits_{k\IN K}\!\!\LRb{z_k{-}x(t_k)}^2\!{+}\\
\!\!{+}\iint\limits_{G}
\!\!\LRs{\alpha_x^2\!\!\LRb{\frac{\partial^2f}{\partial x^2}}^2\!\!\!{+}2\alpha_x\alpha_v\!\!\LRb{\frac{\partial^2f}{\partial x\partial v}}^2 \!\!\!{+}\alpha_v^2\!\!\LRb{\frac{\partial^2f}{\partial v^2}}^2}\!\! dxdv.
\end{array}
\end{equation}
\noindent In this case \SVF-method finds two optimal coefficients $\alpha_x^*$ and $\alpha_v^*$ minimizing the CV-error $\sigma\LRb{\alpha_x,\alpha_v}$, which is calculated similarly to equations (\ref{eq:oscill_xfaKk}). 
Calculating solutions $x^*(t)$, $f^*(x,v)$, $\sigma^*$ and  $\rmse^*$  are similar to those were used for Models 0, 1 (see  equations (\ref{eq:oscill_xfOpt}), (\ref{eq:oscill_xOpt})).
The values $\sigma^*$, $\rmse^*$ are presented in Table \ref{tab:oscill}, isolines of function $f^*(x,v)$ and trajectory $\LRb{x^*(t),\dot{x}^*(t)}$ - in Fig. \ref{fig:oscill_f(x,v)_m2}. The plot of function $x^*(t)$ is skipped, because it is almost the same as in Fig. \ref{fig:oscill_m0_svf} for Model 0. 

\begin{figure}[h]
\centering
\includegraphics[width=0.7\linewidth]{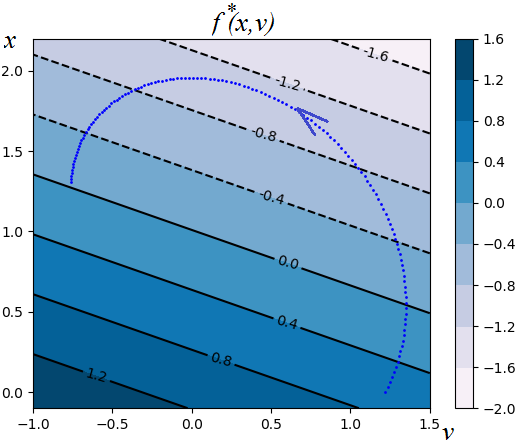}
\caption{Model 2. Isolines of $f^*(x,v)$ and trajectory $\LRb{x^*(t),\dot{x}^*(t)}$.}\label{fig:oscill_f(x,v)_m2}
\end{figure}

Let's draw conclusions from the results obtained above: 
\begin{enumerate}
 \item The errors for Model 2 are almost the same as for Model 0, so, the transition from an arbitrary smooth function to solution of ODE (\ref{eq:oscill_m2}) doesn't reduces capabilities of describing raw data (\ref{eq:oscill_D});
 \item  The level lines in Fig. \ref{fig:oscill_f(x,v)_m2} are parallel and equidistant, so the graph of $f^*(x,v)$ is a plane, 
 so $f^*(x,v)$ is an affine function (a linear one plus a constant).
\end{enumerate}

\paragraph*{Model 3. Oscillator with friction.} 
Affinity of $f(x,v)$ permits to refine right hand of ODE (\ref{eq:oscill_m2}):
\begin{equation}\label{eq:oscill_m3}
\begin{array}{l}
\frac{d^2x}{dt^2}\EQ {-}k\LRb{x-\Dx} - \mu\frac{dx}{dt},
\end{array}
\end{equation}
\noindent where $k,\mu,\Dx$ are unknown parameters to be identified by fitting the raw data (\ref{eq:oscill_D}). 

Results of \SVF-method for Model 3 are presented in Table \ref{tab:oscill}. You can see that the predictive error ($\sigma^*$) is substantially smaller than in other cases, despite that the plot of function $x^*(t)$ is almost equivalent to that of Fig. \ref{fig:oscill_m0_svf}.

Table \ref{tab:oscill} illustrates proposed methodology of quantitative comparing successive model refinements. If the CV-error becomes smaller after the modification, then \quot{we are on the right path}. So, Model 2 is better than Model 0, and Model 3 is better than Model 2. But replacement of Model 0 with Model 1 was wrong.

Root-mean-square error, $\rmse^*$, is another quantitative characteristic of the model - a measure of its \quot{rigidity}. Really, acceptance of \quot{right} assumptions about the model reduces a set of feasible solutions. It leads to refinement of the model (see the change of $\rmse^*$ values in Table \ref{tab:oscill}). Ideally, $\rmse^*$ has to be equal to raw data measurements error. 

Calculations done lead to the conclusion that the Model 3 (oscillator with viscous friction) fits the available data best. Indeed, formula (\ref{eq:oscill}), which was used to generate the data, is one of the solutions of motion equation (\ref{eq:oscill_m3}) - a well-known model of an oscillator with friction.

\begin{table}[h]
\begin{tabular}{clcc}
No & \multicolumn{1}{c}{Model} & $\sigma^*$ & $\mathrm{mse}^*$ \\ \hline
0 & A function. Spline-approximation & 0.1144 & 0.0907 \\ 
1 & 1st Order differential equation & 0.2093 & 0.1869 \\ 
2 & 2nd Order differential equation & 0.1146 & 0.0924 \\ 
3 & Oscillator with friction & 0.1106 & 0.0932 \\ 
\end{tabular}
\caption{Cross-validation error $(\sigma^*)$ and root-mean-square error $(\rmse^*)$ for different models (by oscillator's motion data).}
\label{tab:oscill}
\end{table}

\newcommand{\Ff}{{F_{\scriptscriptstyle F}}}
\newcommand{\Fs}{{F_{\scriptscriptstyle S}}}
\renewcommand{\xaKk}[2]{{x^{#1}_{K{\setminus}#2}}}
\newcommand{\xaK}[2]{{x^{#1}_{#2}}}
\newcommand{\dima}{{n_{\alpha}}} 

\section{Balanced identification method}\label{sec:svf_outline}
In the previous section \SVF-method has been described on model of oscillator with friction. Below a general description of the method is presented. 
It includes formalization of the basic concepts: model; data and model error; main criterion (with regularization) 
and cross-validation scheme.

\vspace{-1em}\paragraph*{Mathematical model} A mathematical model is a set of statements (hypothesis) on the properties 
of the object under study in the form of: equations (e.g. $\dot{x}(t)\EQ f(x)$); 
inequalities (e.g. $\dot{x}(t)\GE0$); statements about the belonging of variables to different spaces (e.g. ${x}(t)\IN C^2[G]$)
(finite dimensional and functional) and sets (e.g. $t\IN T$)); logical expressions (e.g $\exists t:\dot{x}(t)\EQ 0$) etc. 

Let's use a formal description of mathematical model in the form of a system of equations:
\begin{equation}\label{eq:math_model_abc}
\begin{array}{c}
{M(x)=0},
\end{array}
\end{equation}
where each element of the vector $M$ corresponds to one of the model's statements (hypothesis).
For a sake of simplicity the equal sign ($=$) is used in (\ref{eq:math_model_abc}), but it can be substituted by 
various inequalities ($\gtrless$), adhesion sign ($\IN$), 
logical conditions, etc. For example, a model (\ref{eq:oscill_m1}) is written in the form of four statements:

\begin{equation}\label{eq:oscill_m1_abc}
\begin{array}{c} 
\frac{dx}{dt}\EQ f\LRb{x(t)},\\
~x(\cdot)\IN\mathrm{C}^2(T),\\
~f(\cdot)\IN \mathrm{C}^2(X),\\
~t\IN T.
\end{array}
\end{equation}

In many cases the model may be presented in the following form as is customary in optimization literature:
\begin{equation}\label{eq:math_model_vvv}
\begin{array}{c}
x\in Q\eqdef\SET{M(x)=0, ~x\IN S{\subseteq} X},
\end{array}
\end{equation}
\noindent where $x$ is a vector of model variables; $Q$ is a feasible set defined by a set of constraints, where:
$M(\cdot)$ is a system of equations; $X$ is a space of model variables (may be a composition of Euclidean and functional spaces); $S$ is a set defined by additional simple constraints, e.g. inequalities (like ${\lessgtr}0$), logical expressions, integer value condition on variables $x$, etc. 
For instance, the monotonic regression \cite{bib:sysoev2019smoothed} may be considered as special case of (\ref{eq:math_model_vvv}) with inequalities.

The variables $x$ may include both: the functional variables (e.g. the ones describing the trajectory of a dynamic system in Section \ref{sec:oscill}); the model parameters sought (e.g. elastic and viscous friction coefficients in the same demo). 
\vspace{-1em}\paragraph*{Data and model error}

Usually, a set defined in (\ref{eq:math_model_vvv}) is too wide. Identifying the model means selecting the best variables $x\IN Q$ corresponding to the measurements of the object under study. The set of measured data may be described as following:
\begin{equation}\label{eq:mm_D}
\begin{array}{c}
D = \SET{(z_k, P_k(x)): k\IN K}, z_k\IN\RR^{n_z}, P_k{:}X{\to}\RR^{n_z},
\end{array}
\end{equation}
\noindent where: $k$ is the index of the measurement from finite set $K$; $z_k$ is the vector of measured values; $P_k(x)$ is an operator expressing the measured values via model variables. The values of $z_k$ may contain measurements error $\varepsilon_k\EQ\NORM{z_k{-}P_k(x)}$. Below, for simplicity, the source data set $D$ is identified with the set $K$ of indices of pairs $(z_k, P_k(x))$.

\vspace{-1em}\paragraph*{Main criterion}
Consider a problem of minimization of the balanced criterion (model simplicity  vs. data fitting) that depends on $x$, $K$ and a vector of regularization coefficients $\alpha\IN\RR^{\dima}$:
\begin{equation}\label{eq:mm_mainOpt}
\begin{array}{c}
F(x,K,\alpha) = \frac{1}{\MOD{K}}\Ff\LRb{x,K}+\Fs\LRb{\alpha,x}\to\min\limits_{x\IN Q},
\end{array}
\end{equation}
where $\Ff(\cdot)$ is a part of criterion responsible for the data fitting and $\Fs(\cdot)$ is a part that responsible for model complexity (for smoothness as a rule).

Functional $\Ff(x,K)$ is a measure of model approximation error. Usually, additive measure is considered, i.e. if $K\EQ K_1{\bigcup}K_2$, $K_1{\bigcap}K_2\EQ\varnothing$, then $\Ff(x,K)\EQ\Ff(x,K_1){+}\Ff(x,K_2)$ (for any $x$). For definiteness, assume that $\frac{1}{\MOD{K}}\Ff(\cdot)$ is a mean-square error:
\begin{equation}\label{eq:mm_FfSD}
\begin{array}{c}
\Ff(x,K) = \sum\limits_{k\IN K}\LRb{z_k{-}P_k(x)}^2.
\end{array}
\end{equation}

Functional $\Fs(\alpha,x)$ is a measure of model complexity (it is one of possible interpretations of Tikhonov regularization functional \cite{bib:tikhonov1980on}). It depends on the solution sought and the vector of non-negative regularization coefficients $\alpha$. The individual components of $\alpha$ are the penalty coefficients for various aspects of complexity of the solution. 
It is presumed that functional $\Fs(\alpha,x)$ increases monotonically over any $\alpha_i$ (a component of the vector $\alpha$), tends to zero as $\alpha_i{\to}0$, and tends to ${+}\infty$ as $\alpha_i{\to}\infty$. 

The functional of complexity is being chosen  based on the specificity of the object's model (\ref{eq:math_model_abc}). These functionals can be different characteristics of the curvature of functions constituting the model and/or some integral characteristics, such as energy or entropy. For instance, in the demo example in Section \ref{sec:oscill}, for single-variable function $f(x)$ (see Model 1): $\Fs(\alpha,f)\EQ \alpha\int\limits_{X}\!\!\LRb{\frac{d^2f(x)}{dx^2}}^2\!\!dx$, and for two-variable function $f(x,v)$ (see Model 2):
\begin{equation*}
\begin{array}{l}
\!\!\Fs(\alpha,f)\EQ\!\!\!\iint\limits_{X{\times}V}\!\!\LRs{\alpha_x^2\!\!\LRb{\frac{\partial^2f}{\partial x^2}}^2\!\!\!{+}2\alpha_x\alpha_v\!\!\LRb{\frac{\partial^2f}{\partial x\partial v}}^2 \!\!\!{+}\alpha_v^2\!\!\LRb{\frac{\partial^2f}{\partial v^2}}^2}\!\! dxdv.
\end{array}
\end{equation*}

For a fixed $\alpha$ the problem of identification for a given set of measurements $K$ is to find the value $x^{\alpha^*}_K$ that minimizes the functional (\ref{eq:mm_mainOpt}) on a feasible set (\ref{eq:math_model_vvv}):
\begin{equation*}
\begin{array}{l}
\xaK{\alpha}{K}\EQ\Argmin\limits_{x}\SET{F\LRb{x,K,\alpha}: x\IN Q}.\\
\end{array}
\end{equation*}
It is required to find a balanced solution, i.e. such a value $\alpha^*$ that provides a reasonable compromise between the model error and its simplicity.

\vspace{-1em}\paragraph*{Cross-validation procedure}
To find $\alpha^*$, a cross-validation procedure is used \cite{bib:kohavi1995study}. The dataset $K$ is subdivided into disjoint subsets $K_i$ of statistically independent measurements:
\begin{equation}\label{eq:mm_Ki}
\begin{array}{l}
K\EQ \bigcup\limits_{i\IN I}K_i, K_i{\bigcap}K_j\EQ\varnothing, i{\neq}j
\end{array}
\end{equation}
\noindent (let's note that in the case, when measurement errors are random, any non-intersecting sets will be independent).

\newcommand{\SCV}{\sigma(\alpha)}
\newcommand{\SCVa}[1]{\sigma(#1)}
Exclude some subset $K_i$ from the set $K$. Let's find the minimum (\ref{eq:mm_mainOpt}) for a given value $\alpha$ and for remaining set $K{\setminus}K_i$ (training set). Let $\xaKk{\alpha}{K_i}$ be a solution:
\begin{equation}\label{eq:mm_xaKKi}
\begin{array}{l}
\xaKk{\alpha}{K_i}\EQ\Argmin\limits_{x}\SET{F\LRb{x,K{\setminus}K_i,\alpha}: x\IN Q},
\end{array}
\end{equation}
\noindent then the value $\Ff(\xaKk{\alpha}{K_i},K_i)$ is approximation error on test set $K_i$. Repeating this procedure for all subsets $K_i$, $i\IN I$ and summing up the results, we obtain a cross-validation error for a given $\alpha$:
\begin{equation}\label{eq:mm_SCV}
\begin{array}{l}
\SCV=\sum\limits_{i\IN I}\Ff(\xaKk{\alpha}{K_i},K_i).
\end{array}
\end{equation}

The sought vector of optimal weight coefficients and the corresponding solution $x^*$ (for a whole set of measurements $K$) are defined as follows:
\begin{equation}\label{eq:mm_axOpt}
\begin{array}{l}
\alpha^*=\Argmin\limits_{\alpha{\geqq}0}\SCV, \\
x^*=\xaK{\alpha^*}{K}=\Argmin\limits_{x}\SET{F\LRb{x,K,\alpha^*}: x\IN Q}.
\end{array}
\end{equation}

So, the search of optimal coefficients $\alpha$ for a given dataset (\ref{eq:mm_D}) is a \textit{bilevel optimization problem} \cite{bib:dempe2002bilevel}: with $\MOD{I}$ independent optimization problems (\ref{eq:mm_xaKKi}) (these problems may be the variational ones) in the lower level to calculate the value of cross-validation error $\SCV$; with finding in the $\alpha^*$ by minimizing the $\SCV$, (\ref{eq:mm_axOpt}), in the upper level.

In the case when model error $\Ff(\cdot)$ is the standard deviation (\ref{eq:mm_FfSD}) we have \textit{mean (modeling) cross-validation error} ($\sigma^*$) and \textit{root-mean-squared error} ($\rmse^*$):
\begin{equation}\label{eq:mm_CveMse}
\begin{array}{l}
\sigma^*=\sqrt{\FRAC{1}{\MOD{K}}\sum\limits_{i\IN I}\sum\limits_{k\IN K_i}\LRs{z_k-P_k\LRb{\xaKk{\alpha^*}{K_i}}}^2  }, \\
\rmse^*=\sqrt{\FRAC{1}{\MOD{K}}\sum\limits_{k\IN K}\LRs{z_k-P_k(x^*)}^2}.
\end{array}
\end{equation}

The practical applicability of \SVF-method depends on the complexity of the model of the object under study, on the number of measurements and their quality. Generally speaking, \SVF-method may lead to variational problems if the model variables are functions of a continuous argument. In this case, discretization is applied. E.g., continuous domain is replaced with a finite set of points  with unknown values of the sought function ​​at these points. Or unknown function is replaced with polynomials with unknown coefficients. As a result, all problems become finite-dimensional mathematical programming ones. To solve them one can effectively use available solvers and high-performance computing environments.

\newcommand{\an}{{\alpha^{\mathsf{\nu}}}}
\newcommand{\anu}[1]{{\alpha^{#1}}}
\newcommand{\sn}{{\sigma^{\mathsf{\nu}}}}
\newcommand{\snu}[1]{{\sigma^{#1}}}

\section{\SVF-technology as software implementation of \SVF-method }\label{sec:tech}
Method of balanced identification (in short \SVF-method) has been described above in rather abstract manner. 
In this section the main details of its software implementation are presented. Namely the capabilities of symbolic notation of the model's equations and of the method options; discretization of continuous variables, differential and integral equations; surrogate optimization method for solving bilevel optimization problem (\ref{eq:mm_axOpt}); usage of Python Pyomo package and Everest optimization service.

\subsection{\SVF--technology workflow}\label{subsec:svf_wf}
\begin{figure}[h]
\centerline{\includegraphics[scale=0.3]{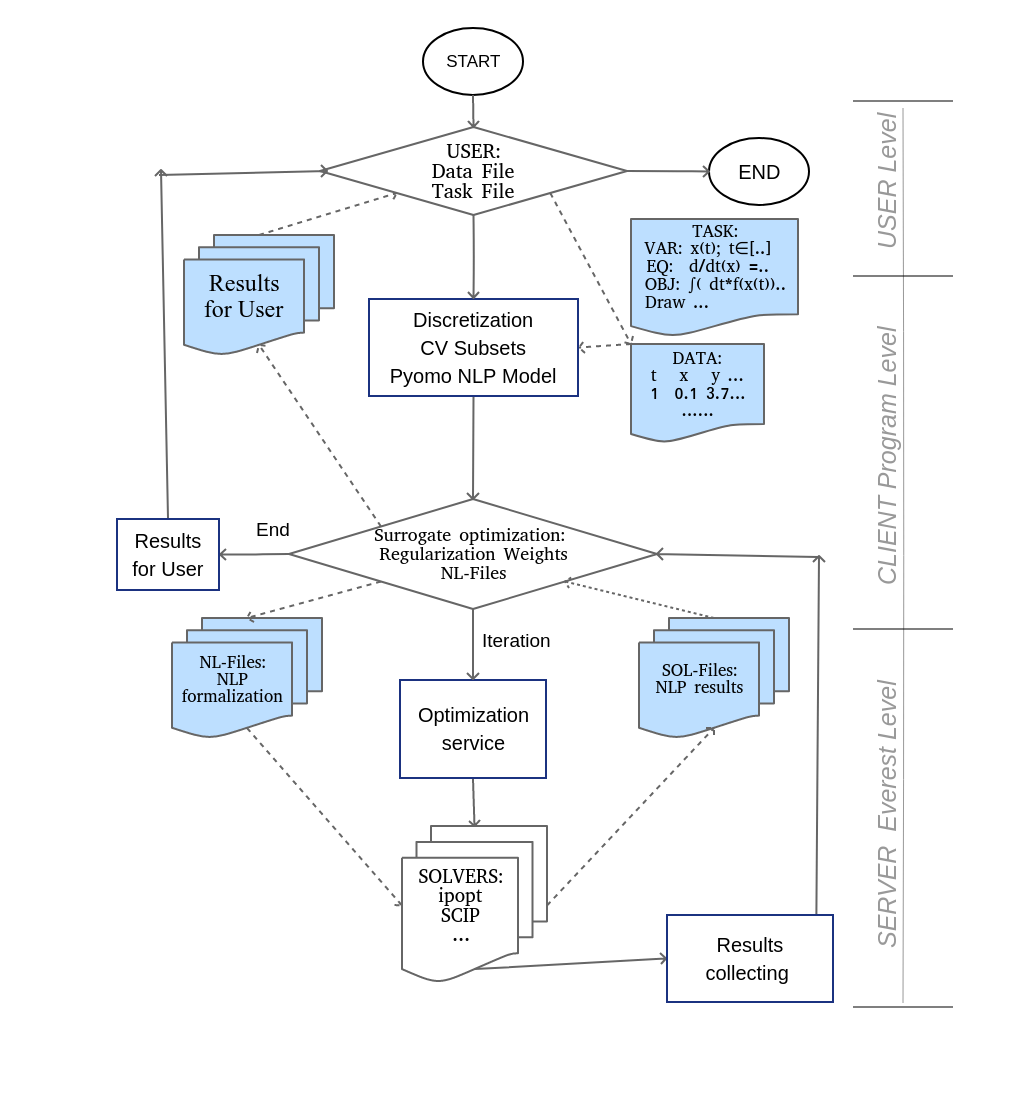}} 
\vspace{-1cm}
\caption{\SVF-technology workflow}
\label{fig:svf_wf}
\end{figure}

The general scheme of the human-computer technology that implements 
the balanced identification method is shown in Fig.~\ref{fig:svf_wf}. 
At the first level (user level), an expert (shape 1 in Fig.~\ref{fig:svf_wf}), who has an idea of ​​the object under study, should prepare a measurement file and a task-file (shape 2 in Fig.~\ref{fig:svf_wf}). 
The data-file contains a table with experimental data (in text format or in MS Excel or MS Acsess formats).

A text task-file usually contains the name of the data-file, a mathematical description of the object (model), 
a list of unknown parameters, cross-validation specifications, etc.
These files are transferred to the client program (shape 3). Here the variational problems are replaced by discretization with NLP (finite dimensional Non-Linear mathematical programming Problem), the sets for CV are generated and the Pyomo NLP Model is formalized. 
The constructed data structures are transmitted to the surrogate 
optimization subroutine (shape 4), which implements an iterative numerical search for unknown model parameters and regularization 
coefficients to minimize the cross-validation error. This subroutine includes parallel solving of mathematical programming problems 
in the distributed environment of Everest optimization services. Pyomo package converts the NLP description in the so-called NL--files 
(shape 5), which at the server level are processed by solvers (shape 7) under the control of the service for solving optimization problems (shape 6). 
The solutions obtained by solvers (shape 8) are collected (shape 9) and sent back to the client level (shape 4), where they are analyzed: 
the conditions for the completion of the iterative process are checked.

If the conditions are not met, surrogate optimization algorithm calculates new values ​​of regularization coefficients and the process repeats. 
Otherwise, the program prepares the results (shape 10), calculates the errors, writes the solution files, 
draws the graphs of the functions found (shape 11), and presents them to the researcher (shape 1).
The obtained results, especially the values of the modeling errors, are used by experts as an argument to choose a new (or modified) model or to stop calculations.

\subsection{Symbolic notation of model and identification problem (task-file)}\label{subsec:mng}

Figure \ref{fig:mng_oscill} shows, as an example, the task-file with the oscillator model (\textbf{Model 2}) description 
and parameters of the numerical method under discussion.

\begin{figure}[h]
\begin{minipage}{\linewidth}
\begin{Verbatim}[frame=single]
CVNumOfIter   20
CVstep        21  
RunSolver     ServerParallel                                
Select x, t  from  ../DATA/Spring5.dat  
GRID:   t  ∈ [ -1.,  2.5, 0.025 ]
        X  ∈ [ -0.1, 2.2, 0.1   ]
        V  ∈ [ -1,   1.5, 0.1   ]
VAR:  x ( t )
      v ( t )
      f ( X, V );  PolyPow = 6
EQ:   d2/dt2(x(t)) == f(x(t),v(t));                        
      v(t) == d/dt(x(t));  
OBJ:  x.MSD() + f.Complexity(Penal[0],Penal[1])
Draw
EOF
\end{Verbatim}
\end{minipage}
\caption{Task-file for oscillator. Model 2}\label{fig:mng_oscill}
\end{figure}

The first three lines define the control parameters: \quot{\texttt{CVNumOfIter}} is the maximum number (20) of iterations (of surrogate optimization procedure), 
\quot{\texttt{CVstep}} is the number (21) of subsets into which the data set is splitted for the cross-validation procedure, 
\quot{\texttt{RunSolver}} specifies the execution mode of calculations (\quot{\texttt{ServerParallel}} means that calculations performed on the server, while independent problems will be solved in parallel). 

The next line defines a data file and names of the columns to be read (“x” and “t”). 
The file contain the dataset defined in (\ref{eq:oscill_D}).

The rest of the task-file contains the \quot{translation} of the mathematical problem (\ref{eq:oscill_m2}),(\ref{eq:oscill_Fxf_m2}) 
into the language of the task-file notation.
In the \quot{\texttt{GRID:}} section the three sets (and a discretization parameters) are specified. 
For example, the expression \quot{\texttt{t∈ [ -1., 2.5, 0.025 ]}} defines a grid from -1. to 2.5 with step 0.025. 

The keyword \quot{\texttt{VAR:}} at the beginning of a line starts a block describing the variables (parameters to be sought).
The grid function \quot{\texttt{x(t)}} is defined on the grid \quot{\texttt{t}} and is linked by the names (\quot{\texttt{x}} and \quot{\texttt{t}}) with the measurements ​​from the file \quot{\texttt{../DATA/Spring5.dat}}.
A function \quot{\texttt{f(X,V)}} is defined as a sixth degree polynomial \quot{\texttt{PolyPow=6}} (with unknown coefficients). 
In this case sets \quot{\texttt{X}}, \quot{\texttt{V}} are used for numerical integration of \quot{complexity} function.  

In the next block (the keyword \quot{\texttt{EQ:}}), two equations (\ref{eq:oscill_m2}) of Model 2 are defined on the grid \quot{\texttt{t}}.

Finally, the key word \quot{\texttt{OBJ:}} defines the objective function consisting of two parts: 
the first one is the data approximation error for the trajectory \quot{\texttt{x(t)}}
(first term in (\ref{eq:oscill_Fxf_m2}), the second one is the measure of the “complexity” 
(second term in (\ref{eq:oscill_Fxf_m2})) of the function {\quot{\texttt{f(X,V)}}}.
The \quot{\texttt{Complexity}} function calculates the regularization penalty for the Penal parameters.
The function can be explicitly specified in the notation.
For example, the first part of the integral looks like 
\quot{\texttt{Penal[0]**2~*∫~d(X)*∫(~d(X)*(d2/dX2(f(X,V)))**2)}}.

The last two lines contain commands to draw the results (functions) and complete the calculations.

Comparing the task-file with the formal mathematical description ({\bf Model 2}), we can note the clarity and convenience of the proposed notation.

\subsection{Implementation of discretization}\label{subsec:discr}
To obtain an approximate solution the continuous (infinite-dimensional, variational) identification problem 
are transformed by discretization into a finite-dimensional one. 
In this balanced identification technology implementation 
the discretization parameters appear (are set) in two objects: 
the discretization step in grids and the polynomial degree for polynomial functions. 
So, continuous functions are replaced by grid or polynomial functions, 
integrals by finite sums, differential equations by sets of algebraic equations.

Relatively simple schemes are used: the trapezoidal integration, the two-point schemes \quot{forward/central/backward} for the first derivatives, the three-point scheme \quot{central} for the second derivative, etc.
The problem is reformulated as a Pyomo model and is written into a special Python-file. A qualified user can always modify it to fit his needs, for example, to use another discretization scheme.
However, there remains the problem of evaluating the proximity of the found numerical solution (based on finite-dimensional mathematical programming problems) to the solution of the original variational problem (in functional space).

\subsection{Implementation by Python \& Pyomo}\label{sec:pyomo}
In the last decade Python became one of the most used programming language in scientific computing in different areas of applied mathematics. The same is true of optimization modeling, i.e. support of AMPL-compatible solvers.

The correct programming of optimization problems (e.g. finite-dimensional analogue of (\ref{eq:mm_xaKKi})) for passing them to a solver via its low-level API may be very time consuming, especially for complex nonlinear problems. The issue is getting more complex when solvers are used in a distributed computing environment. Fortunately, the state-of-the-art solvers
support so called AMPL (A Modeling Language for Mathematical Programming) \cite{bib:ampl-book}. Our research team uses AMPL in our studies on optimization modeling about ten years. In particular, AMPL usage enabled to develop the distributed system of optimization services on the base of Everest toolkit \cite{bib:smirnov2016distributed}. 

Originally the usage of AMPL required commercially licensed translator to generate special NL-file containing all data of the optimization problems to be solved. This NL-file may be passed to AMPL-compatible solver. When solving is successful the solver returns solution as SOL-file, which may be read by AMPL-translator to get the values of variables (including the dual ones). An open source alternative to AMPL-translator appeared in 2012 as Pyomo package \cite{bib:hart2017pyomo}. Pyomo enables to reproduce all \quot{lifecycle} of AMPL-script but without AMPL-translator. Description of optimization problems, generation of NL-files and reading of SOL-files all may be done in any Python application. 

It is well known that Python is not intended for intensive computing. It is rather a scripting language providing data exchange between high-performance software. As to Pyomo, the bottleneck may be in the preparing of NL-file. Fortunately, there are a number of tricks to reduce elapsed time and  a delay become prominent for rather big optimization problems with dozens of thousands of variables and constraints. The most hard computations related to solving of problems presented by NL-files is performed by rather efficient solvers written in C, C++ and Fortran.

One of the results of the task-file processing is a special software module that 
reformulates the optimization problem into Python language using the Pyomo optimization modeling package. 
The fragment of an example of auto-generated Python code corresponding to the first differential equation (\ref{eq:oscill_Fxf_m2}) is shown in Figure \ref{fig:mng_de_qqq}.
\begin{figure}[h]
\begin{minipage}{\linewidth}
\begin{Verbatim}[frame=single]
def EQ0 (Gr,t) :
  return (
    ((x((t+0.025))+x((t-0.025))-2*x(t)) \
           /0.025**2)==f(x(t),v(t))
  )
Gr.conEQ0 = Constraint( \
    myrange(-1.0+0.025,2.5-0.025,0.025), \
    rule=EQ0 )
\end{Verbatim}
\end{minipage}
\caption{First differential equation (\ref{eq:oscill_Fxf_m2}) as a part of Pyomo model}\label{fig:mng_de_qqq}
\end{figure}

\subsection{Surrogate optimization in \SVF-technology}\label{subsec:surrogate}
Bilevel optimization is a well-known hard challenge for researchers in numerical methods \cite{bib:dempe2002bilevel}. Current implementation of \SVF is based on the approach similar to that of \textit{surrogate optimization} \cite{bib:forrester2008engineering}, i.e. upper-level objective function is approximated by its evaluation in a few points. Corresponding set of points is generated successively until some stopping condition is satisfied.

The \SVF-algorithm has two essential features: \\
1) the computational cost of each iteration is rather high (a set of independent NLP problems is to be solved);\\
2) the resulting CV-error estimate may be distorted because of inexact solution of the NLP-problems at low-level. \\
A special surrogate optimization procedure that takes into account these features has been developed. 
The low-level problem is treated as a \quot{black-box} that calculates 
(spending significant computational resources) an approximate result $\sigma$($\alpha$) for the input vector $\alpha$. 

So, we are looking for an approximate minimum of function $\SCV$ (\ref{eq:mm_SCV}) - an optimal value function of the upper-level problem (block 4 in Fig.~\ref{fig:svf_wf}). The value of $\SCV$ depends on solutions of $\MOD{I}$ independent optimization problems (\ref{eq:mm_xaKKi}) for a given $\alpha$ (blocks 5-9 in Fig.~\ref{fig:svf_wf}). The method is based on successive approximation of $\SCV$ by 2nd order polynomial of $\alpha$ on a sequence $\SET{\an}$, $\nu\EQ1,2,\dots$.
The coefficients of this polynomial are recalculated to get more and more better approximation of values $\SCVa{\an}$ for already \quot{passed} points $\an$. Namely, function $\SCV$, $\alpha\IN\RR^\dima$ is approximated by multidimensional polynomial:
\begin{equation}\label{eq:method_Pi}
\hspace{-0.9em}\begin{array}{l}
\Pi(\bm{\pi},\alpha)=\hspace{-1.2em}\sum\limits_{1\LE i\LE j\LE\dima}\hspace{-1.5em}\pi_{ij}\alpha_i\alpha_j + \hspace{-1em}\sum\limits_{1\LE i\LE\dima}\hspace{-1.2em}\pi_{i}\alpha_i + \pi_0,~\mbox{where}\\
\bm{\pi}\EQ\LRb{{\pi_{ij} (1\LE i{<}j\LE\dima)},{\pi_i (i\EQ1{:}\dima)},\pi_0}\IN\RR^{\frac{N(N{+}3)}{2}{+}1},
\end{array}
\end{equation}
\noindent components of $\bm{\pi}$ are changing at runtime of the algorithm.

The algorithm builds sequence of points  $\SET{\an}$, $\nu\EQ1,2,\dots$, stepwise as following. Assume that at step $N$ we have a set of pairs 
$\SET{\LRb{\sigma^{\nu},\an}}_{\nu\EQ1}^{N}$, ($\sigma^{\nu}\EQ\SCVa{\an}$).\footnote{At the initial stage, the new values of the vector $\alpha$ are obtained by small perturbations of all its components}  
In what follows, these pairs will be treated as a set of points in $\RR^{\dima{+}1}$. Without loss of generality (maybe after renumbering) assume that $\snu{N}\EQ\min\limits_{\nu\IN1{:}N}\SCVa{\alpha^\nu}$.

Consider a problem of approximation of these points by a polynomial $\Pi(\bm{\pi},\alpha)$ (\ref{eq:method_Pi}). The smaller the distance between $\an$ and the \quot{best} $\anu{N}$ the more weight the error $\MOD{\sigma^\nu{-}\Pi(\bm{\pi},\an)}$ will have in the following penalty function with regularization term (sum of all squared $\pi_{ij}$ from (\ref{eq:method_Pi}) weighted with coefficient $\mu$):
\begin{equation}\label{eq:method_PiF}
\begin{array}{l}
\sum\limits_{\nu\EQ1}^{N}\FRAC{\LRb{\sn{-}\Pi(\bm{\pi},\an)}^2}{e^{\NORM{\an{-}\anu{N}}}}{+}\mu\LRb{\sum\limits_{i\LE j}\pi_{ij}^2 }\to\min\limits_{\bm{\pi}}.
\end{array}
\end{equation}

This approximation (not interpolation !) does not pass exactly through the points $\SET{\LRb{\sn,\an}}_{\nu\EQ1}^{N}$ but it is close to those points, for which $\sn$ is close to minimum of $\SET{\sn}_{\nu\EQ1}^{N}$. 

It is easy to see that the objective function (\ref{eq:method_PiF}) is strongly convex in variables $\bm{\pi}$ and reaches unique minimum at $\bm{\pi}^*(\mu)$. The algorithm searches $\mu^*$ that minimizes error of $\SCV$ approximation (by polynomial $\Pi(\cdot)$) at the best vector $\anu{N}$. Another bilevel optimization problem arises:
\begin{equation}\label{eq:method_optMu}
\begin{array}{l}
\MOD{\Pi\LRb{\bm{\pi}^*(\mu),\anu{N}} - \snu{N}}\to\min\limits_{\mu\GE0},
\end{array}
\end{equation}
\noindent where $\bm{\pi}^*(\mu)$ is a solution of the lower-level problem (\ref{eq:method_PiF}).

Note, that in the upper-level problem we are looking for an optimal scalar value $\mu^*$ and the lower-level problem is effectively solvable (as the strongly convex one). Thus, optimal solution $\mu^*$ may be found by some of known \textit{line search} algorithms, e.g. \quot{bisection method}.

After polynomial approximation $\Pi\LRb{\bm{\pi}^*(\mu^*),\alpha}$ is determined, the next vector $\alpha^{N{+}1}$ is found as a solution of the following auxiliary problem, similar to the one is used in \textit{Pshenichny-Danilin linearisation method} \cite{bib:pshenichnyi1993linearization}:
\begin{equation}\label{eq:method_nextA}
\begin{array}{l}
\LRb{\nabla_{\alpha}\Pi(\bm{\pi}^*(\mu),\anu{N})}^{\mathsf{T}}\!\!\LRb{\alpha{-}\anu{N}} + \frac{1}{2}\NORM{\alpha{-}\anu{N}}^2,
\end{array}
\end{equation}
\noindent where $\nabla_{\alpha}\Pi(\cdot)$ is a gradient of $\Pi(\bm{\pi},\alpha)$ w.r.t. variable $\alpha$, and $x^{\mathsf{T}}y$ means scalar product of vectors.
Let $\alpha^{N{+}1}$ be a solution of this problem (it exists and is unique). 

At the next step of the algorithm we calculate the value $\snu{N{+}1}\EQ\SCVa{\anu{N{+}1}}$ as a result of solving the set of problems (\ref{eq:mm_xaKKi}). If the difference $\MOD{\snu{N{+}1}{-}\snu{N}}$ is less then some predefined small threshold, then algorithm stops and returns $\anu{N}$ (if $\snu{N{+}1}\GE\snu{N}$) or $\anu{N{+}1}$ (if $\snu{N{+}1}{<}\snu{N}$).
If the above difference in $\SCV$ is higher than a predefined threshold, then the previous step of the algorithm is repeated on the extended set of pairs 
$\SET{\LRb{\sn,\an}}_{\nu\EQ1}^{N{+}1}$. 

\subsection{Solvers and Everest optimization service}\label{subsec:solvers}
As a basic tool for solving independent mathematical programming problems in parallel (shape 6-8 in Fig.~\ref{fig:svf_wf}) we use the Everest optimization service \cite{bib:smirnov2016distributed,bib:smirnov2018dd}. There are a number of services available at the URL \mbox{\url{https://optmod.distcomp.org}}. All of them are based on Everest toolkit \cite{bib:EverestISPDC2015}. Each service, an Everest-application in \quot{Everest-terminology}, is a front-end of a pool of solvers deployed and running on a number of various computing resources (desktops, standalone servers, clusters). This computing environment is extensible because the potential users can add their own computing resources to the optimization environment via the open source Everest agent \mbox{\url{http://everest.distcomp.org/docs}}. 

Current implementation of \SVF-technology is based on \textit{solve-ampl-stub} Everest-application.
This service solves mathematical programming problems 
presented by their NL-files \cite{bib:ampl-book} and returns SOL-files with the solutions found. 
Now the service provides unified access to the following solvers, allowing to solve the main types of mathematical programming problems (LP/MILP/NLP/MINLP):
\begin{itemize}
 \item Ipopt (Coin-OR Interior Point Optimizer, NLP), \mbox{\url{https://github.com/coin-or/Ipopt}};
 \item CBC (Coin-OR Branch-and-Cut, LP, MILP), \mbox{\url{https://github.com/coin-or/Cbc}};
 \item SCIP (Solving Constraint Integer Programs, LP, MILP, MINLP), \mbox{\url{http://scip.zib.de}};
 \item Bonmin (COIN-OR Basic Open-source Nonlinear (convex) Mixed Integer programming, MINLP), \mbox{\url{https://github.com/coin-or/Bonmin}}.
\end{itemize}
\noindent One should say that CBC and Bonmin seems to be a little outdated now. All the results presented above, in Section \ref{sec:oscill}, and below, in Section \ref{sec:use_cases}, are obtained with solver Ipopt. Currently we are beginning to use SCIP global-optimization solver \cite{bib:GleixnerBastubbeEifleretal.2018}. We do it carefully, because the search of global optimum takes much more time than the seacrh of local solution that Ipopt does.

Besides, we are working on implementation of an analogue of AMPLX (an extension of AMPL-scripting to exchange data with Everest optimization services, \cite{bib:smirnov2016distributed}) on the base of Pyomo - so called PyomoEverest, \mbox{\url{https://github.com/distcomp/pyomo-everest}}. It provides developers with simple Python API to exchange data between Pyomo optimization models and Everest optimization services. 
In addition to the \texttt{solve-ampl-stub} 
another Everest-application, \texttt{solve-set-opt-probs}, 
processing a set of independent problems more effectively by the generic Parameter Sweep application \cite{bib:EverestParameterSweep2015} is available now (see the list of applications at the \url{https://optmod.distcomp.org/}).

\section{Successful use cases}\label{sec:use_cases}

The technology has been successfully used in different application areas enlisted below. 

\textbf{Plant physiology} \cite{bib:NSCF2016-Transpiration}. The dataset (about 1600 points in time, step - 30 min) consist of 
  transpiration and climatic (photoactive radiation, humidity and temperature) measurements. 
  Among a number of (about 10) models of different complexity, a dynamic model of the water in the plant was chosen.
  
\textbf{Ecology} \cite{bib:sokolov2019application}. The method of balanced identification was used to describe 
  the response of Net Ecosystem Exchange of $\mathrm{CO}_2$ (NEE) 
  to the change of environmental factors, and to fill the gaps in continuous $\mathrm{CO}_2$ flux measurements in a sphagnum peat bog in the Tver region.
  The dataset is about 2000 points (step - 30 min). 
  The selected model (the dependence of NEE on climatic parameters) was a superposition of one- and two-variable functions.

\textbf{Meteorology} \cite{bib:lavruhin2016radon}. The model of radon dynamics in the atmosphere was build. 
  The dataset (about 3000 points, about month, step 15 min) contains measurements of characteristics of the $\gamma$-radiation field by a gamma spectrometer.
  The model, consisting of 2 ODEs, describes the dynamics of the volume activity of ${}^{222}\mathrm{Rn}$ 
  and its daughter products ${}^{214}\mathrm{Pb}$ and ${}^{214}\mathrm{Bi}$.

\textbf{Pollution, geography, soil erosion} \cite{bib:linnik-2016-geokhi}. 
  The datasets (from 100 to 1000 points by X and Y coordinates, step - about 100 meters), 
  consist of altitude (SRTM satellite data) and ${}^{137}\mathrm{Cs}$ decay activity (airborne gamma data) measurements. 
  A number of models were considered (functional dependencies), including digital terrain, soil erosion and ${}^{137}\mathrm{Cs}$ contamination models. 

\textbf{Air pollution} \cite{bib:EGU2019-AirPolution}.
  The problem is considered of selection of mathematical description based on 2D diffusion-transport model 
  for the solution of the inverse problem of atmospheric pollution dispersion.
  The following problems were studied (see poster at \url{https://www.researchgate.net/publication/334523588}):

    F)  The estimation of distribution and sources of the atmospheric aerosol pollution (PM10) in the North of France. 
    The dataset (about 30000 point) consists of 12 stations of Atmo Hauts-de-France pollution measurement network.
    
    R)  The evaluation of pollutant emissions and the average wind rose for contamination sources with known positions. 
    Measurements of pollution of soil and lake’s water by metals (Cu, Ni, etc.) were obtained 
    in the Kola Peninsula in Russia at about 100 stations.
  
\textbf{Medicine} \cite{bib:CMMASS2019-Glucose}. Different models of glucose and insulin dynamics in human blood are analyzed.
  The datasets (about 20 measurements of glucose concentration and insulin) contain intra venous glucose tolerance test results. 

\textbf{Thermonuclear plasma} \cite{bib:afanasiev2017inverse,bib:VANT-2017-Plasma}. 
  We compared several models (integro-differential, partial derivatives) describing the initial stage of the fast nonlocal transport events, 
  which exhibit immediate response, in the diffusion time scale, 
  of the spatial profile of electron temperature to its local perturbation, 
  while the net heat flux is directed opposite to ordinary diffusion (i.e. along the temperature gradient). 
  The results of experiments on a few thermonuclear installations were used. The datasets contain about 100 points.
  
\textbf{Evolutionary theory} \cite{bib:sokolov2016evolution}. 
  The optimization problem, which reflects the mechanisms of regulation of the speed of evolution 
  that provide an adequate population response to the direction and rate of environmental change, was considered. 
  The numerical experiment results show plausible dependences of age-specific fertility on the rate of environmental changes and also describe and explain a number of evolutionary effects.

\section{Discussion}\label{sec:concl}
Presented \SVF-method and its software implementation are an interesting combination of knowledge from the areas of  data analysis, optimization and distributed computing software including: cross-validation and regularization methods, algebraic modeling in optimization and methods of optimization, automatic discretization of differential and integral equation, optimization REST-services, etc. 

Current implementation of balanced identification method is based on state-of-the-art open source NLP-solvers (AMPL-compatible): Ipopt \cite{bib:wachter2006ipopt} and SCIP \cite{bib:GleixnerBastubbeEifleretal.2018}. Optionally, optimization services \cite{bib:smirnov2016distributed} of extensible computing environment \mbox{\url{https://optmod.distcomp.org}} based on Everest toolkit \cite{bib:EverestISPDC2015} may be used as well.

The implementation supports symbolic notation to describe the model and corresponding identification problem. This feature automates the formulation of optimization problems with less user effort. A special module generates fragments of Python codes based on a special discretization subsystem and Pyomo framework to interact with open source AMPL-compatible solvers. 

Success story list of the presented approach already is rather long. Nevertheless, authors realize its domain of applicability and important unresolved issues. 

First of all, the method does not provide a fully automated \quot{model generation and/or selection}. \SVF-technology enables quantitative comparison within a given set of models, but it is the researcher who make modifications of the model (e.g. change model constraint equations).

The second issue concerns the search of global optimum in problems (\ref{eq:mm_xaKKi}), which is required for cross-validation. In general case, these problems may have a lot of local extrema, moreover, they may even contain integer variables. If so, NLP-solver Ipopt can not guaranty a global optimum. It is possible to use SCIP solver implementing branch-and-bound algorithm for mixed-integer nonlinear mathematical programming problems, but at the cost of much more longer solving time (and more memory consumption). A parallel extensions of SCIP are available \cite{bib:shinano2011parascip,bib:shinano2017fiberscip}, but first we are going to use simpler workarounds (e.g. Ipopt multi-start options and some easily verifiable necessary conditions of that Ipopt returns real global optimum). In current practice some \quot{informal} analysis of \SVF-algorithm convergence with respect to decreasing of modeling cross-validation error is applied. In the nearest future we plan to perform computational experiment where global optimization solver SCIP and its parallel extensions FiberSCIP or ParaSCIP will be used alongside with Ipopt. E.g. one call of SCIP per a ten calls of Ipopt.

\newcommand{\Nx}{{N_x}}
\newcommand{\Ny}{{N_y}}
Discretization is the third important challenge for \SVF-method. Presence of ordinary and partially differential equations, along with integral ones, is typical for models in physics, biology, medicine etc. Current implementation is based on rather simple discretization module that implements two basic schemes: finite-difference and approximation by sums of polynomials with unknown coefficients. We plan to analyse capabilities of DAE subsystem of Pyomo (discretization of Differential and Algebraic Equations) for \SVF-implementation, \cite{bib:nicholson2018pyomo}. One should remind that \SVF implementation already uses Pyomo package for programming of optimization problems and for data exchange with AMPL-compatible solvers.

One of the issue related to discretization concerns the case when the model includes composition $f{\circ}x$ of unknown functions $f(y)$, $x(t)$. E.g., the Model 1 (\ref{eq:oscill_m1}) includes differential equation $\dot{x}{=}f(x(t))$ with unknown functions $f(y)$, $x(t)$. Current implementation of \SVF uses the following approach: inner function $x(t)$ is approximated by a grid function and outer $f(y)$ - by polynomial:
\begin{equation}\label{eq:fPoly_xGrid}
\begin{array}{l}
x(t)\sim \SET{\LRb{x_i,t_i}: i=0{:}\Nx},\mbox{~where~} x(t_i)=x_i;\\
f(y)\sim \sum_{k=0}^{K}{{p}_{k}}\cdot {{y}^{k}}.
\end{array}
\end{equation}
\noindent Here unknown variables are $\SET{x_i{:}i\EQ0{:}\Nx}$ and $\SET{p_k,k\EQ1{:}K}$. Discretization of differential equation (\ref{eq:oscill_m1}) looks like a number of equations:
\begin{equation}
\frac{x(t_i){-}x(t_{i{-}1})}{\Delta_i t}=\sum\limits_{k=0}^{K}{{p}_{k}}\cdot{x(t_i)^{k}}.
\end{equation}
We plan to use one more approach based on a smoothed piecewise approximation (e.g., see \cite{bib:Chen1996,bib:Zhou2019}) of outer function approximated by a grid as well.

Another issue concerns the domain of applicability of discretization itself, e.g. for the models based on optimal control problems. Here the applicability of discretization depends on whether the corresponding problem has any special properties, like a \quot{chattering phenomenon} of a global optimal solution. 
The well-known Fuller problem \cite{bib:fuller1960relay} is one of the simplest demonstration of this effect. At first glance, it is a very simple optimal control problem on a finite time interval, but for some segment of initial conditions any optimal solution can not be well approximated by any finite discretization of time interval.


\section*{Acknowledgments}
The results of section 4 were obtained within the Russian Science Foundation grant (project No. \text{16-11-10352}). The results of sections 1-3, 5 were obtained within Russian Foundation for Basic Research grant (project No. \text{20-07-00701}).

\bibliographystyle{unsrt}
\bibliography{svf}

\begin{thebibliography}{10}

\bibitem{bib:sokolov2018choice}
A.V. Sokolov and V.V. Voloshinov.
\newblock Choice of mathematical model: balance between complexity and
  proximity to measurements.
\newblock {\em International Journal of Open Information Technologies}, 6(9),
  2018.

\bibitem{bib:tikhonov1980on}
A.N. Tikhonov.
\newblock On mathematical methods for automating the processing of
  observations.
\newblock In {\em Problems of Computational Mathematics}, pages 3--17, 1980.

\bibitem{bib:kohavi1995study}
Ron Kohavi et~al.
\newblock A study of cross-validation and bootstrap for accuracy estimation and
  model selection.
\newblock In {\em Ijcai}, volume~14, pages 1137--1145. Montreal, Canada, 1995.

\bibitem{bib:hastie2009elements}
Trevor Hastie, Robert Tibshirani, and Jerome Friedman.
\newblock {\em The elements of statistical learning: data mining, inference and
  prediction}.
\newblock Springer, 2 edition, 2009.

\bibitem{bib:kuhn2013applied}
Max Kuhn and Kjell Johnson.
\newblock {\em Applied predictive modeling}, volume~26.
\newblock Springer, 2013.

\bibitem{bib:rozhenko2005theory}
A.I. Rozhenko.
\newblock {\em {Theory} and {Algorithms} of {Variational}
  {Spline}-{Approximations}}.
\newblock Novosibirsk State Technical University, 2005.
\newblock (in Russian).

\bibitem{bib:hardle1990applied}
Wolfgang H{\"a}rdle.
\newblock {\em Applied nonparametric regression}.
\newblock Number~19. Cambridge university press, 1990.

\bibitem{bib:strijov2010nonlinear}
Vadim Strijov and Gerhard~Wilhelm Weber.
\newblock Nonlinear regression model generation using hyperparameter
  optimization.
\newblock {\em Computers \& Mathematics with Applications}, 60(4):981--988,
  2010.

\bibitem{bib:sysoev2019smoothed}
Oleg Sysoev and Oleg Burdakov.
\newblock A smoothed monotonic regression via {L2} regularization.
\newblock {\em Knowledge and Information Systems}, 59(1):197--218, 2019.

\bibitem{bib:dempe2002bilevel}
Stephan Dempe.
\newblock {\em Foundations of bilevel programming}.
\newblock Springer Science \& Business Media, 2002.

\bibitem{bib:pshenichnyi1993linearization}
B.N. Pshenichnyi and A.A. Sosnovsky.
\newblock The linearization method: Principal concepts and perspective
  directions.
\newblock {\em Journal of Global Optimization}, 3(4):483--500, 1993.

\bibitem{bib:boyd2011distributed}
S.~Boyd, N.~Parikh, E.~Chu, B.~Peleato, and J.~Eckstein.
\newblock Distributed optimization and statistical learning via the alternating
  direction method of multipliers.
\newblock {\em Foundations and Trends in Machine Learning}, 3(1):1--122, 2011.

\bibitem{bib:smirnov2016distributed}
S.~Smirnov, V.~Voloshinov, and O.~Sukhoroslov.
\newblock Distributed optimization on the base of {AMPL} modeling language and
  {Everest} platform.
\newblock {\em Procedia Computer Science}, 101:313--322, 2016.

\bibitem{bib:smirnov2018dd}
S.~Smirnov and V.~Voloshinov.
\newblock On domain decomposition strategies to parallelize branch-and-bound
  method for global optimization in {Everest} distributed environment.
\newblock {\em Procedia Computer Science}, 136:128--135, 2018.

\bibitem{bib:EverestISPDC2015}
O.~Sukhoroslov, S.~Volkov, and A.~Afanasiev.
\newblock A web-based platform for publication and distributed execution of
  computing applications.
\newblock In {\em Parallel and Distributed Computing (ISPDC), 2015 14th
  International Symposium on}, pages 175--184, June 2015.

\bibitem{bib:ampl-book}
R.~Fourer, D.M. Gay, and B.W. Kernighan.
\newblock {\em {AMPL}: A Modeling Language for Mathematical Programming. Second
  edition}.
\newblock Duxbury Press/Brooks/Cole Publishing Company, 2003.
\newblock \url{https://ampl.com/resources/the-ampl-book}.

\bibitem{bib:hart2017pyomo}
W.E. Hart, C.D. Laird, J.P. Watson, D.L. Woodruff, G.A. Hackebeil, B.L.
  Nicholson, and J.D. Siirola.
\newblock {\em Pyomo--optimization modeling in {Python}. 2nd edition},
  volume~67.
\newblock Springer, 2017.

\bibitem{bib:forrester2008engineering}
Alexander Forrester, Andras Sobester, and Andy Keane.
\newblock {\em Engineering design via surrogate modelling: a practical guide}.
\newblock John Wiley \& Sons, 2008.

\bibitem{bib:GleixnerBastubbeEifleretal.2018}
A.~Gleixner, M.~Bastubbe, L.~Eifler, T.~Gally, G.~Gamrath, R.~L. Gottwald,
  G.~Hendel, C.~Hojny, T.~Koch, M.~E. L{\"u}bbecke, S.~J. Maher,
  M.~Miltenberger, et~al.
\newblock The {SCIP} {O}ptimization {S}uite 6.0.
\newblock Technical Report 18-26, ZIB, Takustr. 7, 14195 Berlin, 2018.

\bibitem{bib:EverestParameterSweep2015}
S.~Volkov and O.~Sukhoroslov.
\newblock A generic web service for running parameter sweep experiments in
  distributed computing environment.
\newblock {\em Procedia Computer Science}, 66:477--486, 2015.

\bibitem{bib:NSCF2016-Transpiration}
A.V. Sokolov, V.K. Bolondinsky, and V.V. Voloshinov.
\newblock Technologies for constructing mathematical models from experimental
  data: applying the method of balanced identification using the example of
  choosing a pine transpiration model.
\newblock In {\em National Supercomputer Fjrum (NSCF-2018)}, 2018.

\bibitem{bib:sokolov2019application}
A.V. Sokolov, V.V. Mamkin, V.K. Avilov, D.L. Tarasov, Y.A. Kurbatova, and A.~V.
  Olchev.
\newblock Application of a balanced identification method for gap-filling in
  {CO}$_2$ flux data in a sphagnum peat bog.
\newblock {\em Computer Research and Modeling}, 11(1):153--171, 2019.

\bibitem{bib:lavruhin2016radon}
Yu.E. Lavruhin, A.V. Sokolov, and D.S. Grozdov.
\newblock Monitoring of volume activity in the atmospheric surface layer based
  on the testimony of the spectrometer seg-017: error analysis.
\newblock In {\em Radioactivity after nuclear explosions and accidents:
  consequences and ways to overcome}, pages 359--368, 2016.

\bibitem{bib:linnik-2016-geokhi}
V.G. Linnik, A.V. Sokolov, and I.V. Mironenko.
\newblock 137cs patterns and their transformation in landscapes of the opolye
  of the bryansk region.
\newblock {\em Modern trends in the development of biogeochemistry}, pages
  423--434, 2016.

\bibitem{bib:EGU2019-AirPolution}
A.V. Sokolov, A.A. Sokolov, and Hervé Delbarre.
\newblock Method of balanced identification in the inverse problem of transport
  and diffusion of atmospheric pollution.
\newblock In {\em EGU2019-15175}, volume~26, 2019.

\bibitem{bib:CMMASS2019-Glucose}
A.V. Sokolov and L.A. Sokolova.
\newblock Building mathematical models: quantifying the significance of
  accepted hypotheses and used data.
\newblock In {\em XXI International Conference on Computational Mechanics and
  Modern Applied Software Systems (CMMASS’2019)}, pages 114--115, 2019.

\bibitem{bib:afanasiev2017inverse}
A.P. Afanasiev, V.V. Voloshinov, and A.V. Sokolov.
\newblock Inverse problem in the modeling on the basis of regularization and
  distributed computing in the {Everest} environment.
\newblock In {\em CEUR Workshop Proceedings}, pages 100--108, 2017.

\bibitem{bib:VANT-2017-Plasma}
A.B. Kukushkin, A.A. Kulichenko, P.A. Sdvizhenskii, A.V. Sokolov, and V.V.
  Voloshinov.
\newblock A model of recovering parameters of fast non-local heat transport in
  magnetic fusion plasma.
\newblock {\em Problems of Atomic Science and Technology, Ser. Thermonuclear
  Fusion}, 40(1):45--55, 2017.

\bibitem{bib:sokolov2016evolution}
A.V. Sokolov.
\newblock Mechanisms of regulation of the speed of evolution: The population
  level.
\newblock {\em Biophysics}, 61(3):513--520, 2016.

\bibitem{bib:wachter2006ipopt}
A.~W{\"a}chter and L.T. Biegler.
\newblock On the implementation of an interior-point filter line-search
  algorithm for large-scale nonlinear programming.
\newblock {\em Mathematical programming}, 106(1):25--57, 2006.

\bibitem{bib:shinano2011parascip}
Y.~Shinano, T.~Achterberg, T.~Berthold, S.~Heinz, and T.~Koch.
\newblock {ParaSCIP}: a parallel extension of {SCIP}.
\newblock In {\em Competence in High Performance Computing 2010}, pages
  135--148. Springer, 2011.

\bibitem{bib:shinano2017fiberscip}
Yuji Shinano, Stefan Heinz, Stefan Vigerske, and Michael Winkler.
\newblock {FiberSCIP} -- a shared memory parallelization of {SCIP}.
\newblock {\em INFORMS Journal on Computing}, 30(1):11--30, 2017.

\bibitem{bib:nicholson2018pyomo}
B.~Nicholson, J.D. Siirola, J.-P. Watson, V.M. Zavala, and L.T. Biegler.
\newblock pyomo.dae: a modeling and automatic discretization framework for
  optimization with differential and algebraic equations.
\newblock {\em Mathematical Programming Computation}, 10(2):187--223, 2018.

\bibitem{bib:Chen1996}
Chunhui Chen and O.~L. Mangasarian.
\newblock A class of smoothing functions for nonlinear and mixed
  complementarity problems.
\newblock {\em Computational Optimization and Applications}, 5(2):97--138,
  1996.

\bibitem{bib:Zhou2019}
Zhengyong Zhou and Yunchan Peng.
\newblock The locally {Chen}--{Harker}--{Kanzow}--{Smale} smoothing functions
  for mixed complementarity problems.
\newblock {\em Journal of Global Optimization}, 74(1):169--193, 2019.

\bibitem{bib:fuller1960relay}
A.T. Fuller.
\newblock Relay control systems optimized for various performance criteria.
\newblock volume~1, pages 520--529. Elsevier, 1960.
\newblock \url{https://doi.org/10.1016/S1474-6670(17)70097-3}.

\end{thebibliography}

\end{document}